\documentclass[reqno]{amsart}
\usepackage[margin = 1.5in]{geometry}
\usepackage{amsmath, amssymb, amsthm, fancyhdr, verbatim, graphicx, stmaryrd}
\usepackage{enumerate}
\usepackage[all]{xy}
\usepackage[usenames,dvipsnames]{xcolor}
\usepackage{mathrsfs}
\usepackage{tikz-cd}
\usepackage{framed, hyperref}
\usepackage[titletoc]{appendix}
\usepackage{bbm}
\usepackage{lipsum}
\usepackage{adjustbox}
\usepackage{devanagari}

\numberwithin{equation}{section}

\numberwithin{equation}{section}

\DeclareFontFamily{U}{matha}{\hyphenchar\font45}
\DeclareFontShape{U}{matha}{m}{n}{
      <5> <6> <7> <8> <9> <10> gen * matha
      <10.95> matha10 <12> <14.4> <17.28> <20.74> <24.88> matha12
      }{}
\DeclareSymbolFont{matha}{U}{matha}{m}{n}
\DeclareFontFamily{U}{mathx}{\hyphenchar\font45}
\DeclareFontShape{U}{mathx}{m}{n}{
      <5> <6> <7> <8> <9> <10>
      <10.95> <12> <14.4> <17.28> <20.74> <24.88>
      mathx10
      }{}
\DeclareSymbolFont{mathx}{U}{mathx}{m}{n}

\DeclareMathSymbol{\obot}         {2}{matha}{"6B}

\RequirePackage{xspace}

\usepackage{color}

\newcommand{\F}{\mathbf{F}}
\newcommand{\RR}{\mathbf{R}}

\newcommand{\G}{\mathbf{G}}

\newcommand{\QQ}{\mathbb{Q}}
\newcommand{\Z}{\mathbf{Z}}
\newcommand{\ZZ}{\mathbb{Z}}

\newcommand{\mf}[1]{\mathfrak{#1}}

\newcommand{\Cl}{\operatorname{Cl}}

\newcommand{\Gal}{\operatorname{Gal}}

\newcommand{\R}{\mathbb{R}}

\newcommand{\ol}[1]{\overline{#1}}
\newcommand{\wh}[1]{\widehat{#1}}

\newcommand{\Cal}[1]{\mathcal{#1}}
\newcommand{\A}{\mathbf{A}}

\newcommand{\co}{\colon}

\newcommand{\msf}[1]{\mathsf{#1}}
\newcommand{\bs}{\backslash}

\newcommand{\cHom}{\Cal{H}om}

\DeclareMathOperator{\GL}{GL}

\DeclareMathOperator{\ab}{ab}

\DeclareMathOperator{\Spec}{Spec\,}

\DeclareMathOperator{\Pic}{Pic}

\DeclareMathOperator{\Sym}{Sym}

\DeclareMathOperator{\Ch}{Ch}
\DeclareMathOperator{\Fun}{Fun}

\newcommand{\limit}{\varprojlim}
\newcommand{\colim}{\varinjlim}

\newcommand{\bu}{\bullet}

\newcommand{\dab}{{ab,\bullet}}

\DeclareMathOperator{\op}{op}

\renewcommand\a\alpha
\renewcommand\b\beta
\newcommand\g\gamma
\renewcommand\d\delta
\newcommand\D\Delta

\newcommand\cF{\mathcal{F}}

\newcommand\cO{\mathcal{O}}

\newcommand{\isoarrow}{{~\overset\sim\longrightarrow~}}

\newcommand\ra{\rightarrow}

\newcommand\rC{\mathrm{C}}

\newcommand\rH{\mathrm{H}}

\newtheorem{thm}{Theorem}[section]

\newtheorem{defn-prop}[thm]{Definition-Proposition}

\theoremstyle{remark}

\newtheorem{remark}[thm]{Remark} 

\newtheorem{const}[thm]{Construction}

\def\vv{\vskip10pt}
\makeatletter
\def\th@remark{%
  \thm@headfont{\bfseries}%
  \normalfont 
  \thm@preskip \thm@preskip 
  \thm@postskip\thm@preskip
}
\def\imod#1{\allowbreak\mkern5mu({\operator@font mod}\,\,#1)}
\makeatother

\numberwithin{equation}{section}

\title[Derived class field theory]{Derived class field theory}

\author{Tony Feng}
\thanks{T. Feng was partially supported by an NSF Postdoctoral Fellowship under Grant No. DMS1902927.
}

\address{Tony Feng\\Dept. of Mathematics, University of California at Berkeley,
Berkeley, CA 94720-
}
 \email{fengt@berkeley.edu}

\author{Michael Harris}
\thanks{M. Harris  was partially supported by NSF Grant DMS-2001369, and by a Simons Foundation Fellowship, Award number  663678}

\address{Michael Harris\\
Department of Mathematics, Columbia University, New York, NY  10027, USA}
 \email{harris@math.columbia.edu}

 \author{Barry Mazur}
\thanks{B. Mazur was partially supported by the NSF Grant DMS 2152149}

\address{Barry Mazur\\
Department of Mathematics, Harvard University, Cambridge, MA  02138, USA}
 \email{mazur@g.harvard.edu}

\begin{document}

\maketitle


   \tableofcontents
   

   \section{Introduction}

Our article  is in memory of John Coates, in memory of his energy, his generosity of thought,  his appreciation of ideas.\vv
                                            \centerline{({\it Barry M.:})}

\begin{quote}  He was an inspiration to me from the earliest days that I knew him---when---beginning  in 1969--- he was a Benjamin Pierce Assistant Professor at Harvard, to the later years,  when he was based in Orsay, France and I was at the IHES and when the two of us would jog together as  he would explain his latest mathematical thought.\end{quote} \vv
   \centerline{({\it Michael H.:})}
   
\begin{quote} I had the privilege of meeting John as a new Ph.D.   My thesis had been largely inspired by the Coates-Wiles paper, which had appeared just one year before.  John's  encouragement and support was precious at this stage of my career.  Although I failed to find interesting applications of my thesis work and shifted my attention to other questions, I remained in close contact with John, especially after John moved to France.  Always elegant, always diplomatic, always with just the slightest trace of a smile on his lips, during his brief stay at Orsay and the  {\'E}cole Normale Sup{\'e}rieure John left an influence on  number theory in France that is still felt today.  After he settled in Cambridge he did the same for Europe as a whole.  John was uniquely effective in helping to build European number theory, and this left a deep impression on me when he moved to France, showing me that it was possible to use the modest powers of a European academic creatively as well as constructively.  I never made any  decisions that might significantly affect our mathematical community without first consulting John.
  Inevitably, John's influence on me was primarily mathematical, through his own work and through that of his mathematical descendents.  No other number theorist of his generation had such a vast mathematical family as John---I have published papers with eight of them, with more on the way.   Inspired by Barry's work on Selmer groups and by unpublished work of Ralph Greenberg, John, together with his student Bernadette Perrin-Riou, reformulated and reinvigorated (classical) Iwasawa theory by  extending it to motives.   This perspective shaped my  return to Iwasawa theory, starting with a joint paper with Jacques Tilouine -- the last of John's students in France.  And, though this may not be so immediately apparent, it also shapes my thoughts about the project on which this paper is a report.\end{quote}\vv
  
       \centerline{({\it Tony F.:})}
     \begin{quote}  Unfortunately I never had the pleasure of meeting John Coates in person, but I have had many encounters with his mathematics, which was and continues to be an inspiration for me.  \end{quote}

\section{Beyond class field theory}
Let $F$ be a number field and $K$ be an open subgroup of the restricted product $\prod'_{\mf{p} \in |F|} \cO_{\mf{p}}^{\times}$. Let $K_\infty$ denote the maximal compact subgroup of $(F\otimes_\QQ \RR)^\times$; thus $K_\infty \isoarrow (\pm 1)^{r_1} \times (S^1)^{r_2}$, where $S^1$ is the unit circle; here as usual, $r_1$ is the number of real embeddings and $r_2$ is the number of complex embeddings.
Define the \emph{idele class groupoid} as the quotient stack 
\[
I_K = I_{F,K} := [F^{\times} \bs \A_F^{\times} /K\cdot K_\infty].
\]
Here the brackets mean that we take the quotient in the sense of groupoids, or in other words we form the \emph{homotopy quotient}. The ``idele class group'' traditionally considered in class field theory is the quotient group $F^{\times}  \bs \A_F^{\times} /K\cdot K_\infty$, which can be thought of as $\pi_0(I_K)$. However, the groupoid $I_K$ has an interesting homotopy type that we will also want to consider. 

Traditional class field theory describes  $\pi_0(I_K)$, as $K$ varies, as the abelianizations of certain Galois groups of the maximal field extension of $F$ with ramification dictated by $K$. However, we have good reason to want to describe the entirety of $I_K$, and not just its component groups, in terms of Galois theory. In other words, we would like to enlarge class field theory to account for the entire id\`ele class group and not just its group $\pi_0(I_K)$ of connected components.  The purpose of the present article is to explain that this is possible: what we shall see is that within this new framework, $I_{F,K}$ accounts completely for what we call the {\it derived abelianization} of the absolute Galois group of $F$.

Why might we want this? $I_K$ is  the {locally symmetric space} associated to the reductive group $G = \GL(1)$ over $F$. We expect--- thanks to ideas of Galatius-Venkatesh \cite{GV} --- that there's an  analogous but perhaps  subtler relation between the topology of  locally symmetric spaces  attached to more general reductive groups and  their corresponding  Galois  representations. See \cite{FH} for an introduction to this circle of ideas. The general conjectures seem intractable at present, but they encompass the case $G = \GL(1)$, so we might as well try to solve that (easiest) case first.  This will be done in forthcoming joint work of the authors and Arpon Raksit \cite{FHMR}. Although the case of $\GL(1)$ is relatively simple, it arises as a useful technical tool in studying other $G$ (for example, one might want to ``twist by characters'', or ``fix determinants'', etc), and so nailing down this case should help in the more interesting cases as well. The present survey explains a component of \cite{FHMR} that we call ``derived class field theory''.

In this survey we aim to give an informal and intuitive explanation, therefore omitting technicalities on higher category theory and homotopical algebra, as well as focusing on simplifying special cases. Precise and complete details will appear in the article (joint with Raksit) \cite{FHMR}. 

\section{A derived Langlands correspondence for $\GL(1)$  -- the Galois side }\label{derGalois}







\subsection{Derived abelianization}  Let 
$\Gamma$ be a discrete  group. The abelianization of  $\Gamma$ is an abelian group $\Gamma^{\ab}$  for which the projection $G\ \longrightarrow\ G^{\ab}$  is the universal solution to the problem of morphisms from $G$ to any  abelian group.  That is, for an abelian group $A$,  $$Hom_{\rm gps}(G, A) = Hom_{\rm gps}(G^{\ab}, A).$$ We have \begin{equation}\label{h1}G^{\ab} = G/[G,G] = H_1(G, {\Z}).\end{equation}

Just as $G^{\ab}$ gives us  one-dimensional homology of $G$  (as in Equation (\ref{h1}) above)  the   \emph{derived abelianization} of  $\Gamma$, denoted $\Gamma^\dab$, is  represented by a simplicial abelian group  that is constructed canonically in the appropriate category and captures all of $H_*(G,{\Z})$.    Specifically, there is a canonical isomorphism \begin{equation}\label{pH}\pi_i (\Gamma^\dab) \simeq H_{i+1}(G, {\Z})\end{equation} for $\i\ge 0$.


Intuitively speaking, the derived abelianization should be a kind of ``derived functor of abelianization''. However, the process of ``deriving'' the abelianization functor cannot be approached as in classical homological algebra, since the category of groups is far from being the sort of abelian category to which the classical theory of derived functors applies. What one uses instead is Quillen's theory of \emph{homotopical algebra} \cite{Q}.

Recall that homological algebra is implemented using the notion of chain complex, which however is very specific to abelian categories. In contrast, Quillen's homotopical algebra uses the notion of \emph{simplicial objects}, which applies to totally general categories. A simplicial object of a category $\msf{C}$ is a collection of objects $C_n \in \msf{C}$ for $n \geq 0$, together with maps $C_m \rightarrow C_n$ whose combinatorics is modeled on the maps of standard simplices
\[
\Delta_n := \{(t_0, \ldots, t_n) \in \R^{n+1}_+ \co \sum t_i = 1 \}.
\]
More precisely, the simplex $\Delta_n$ is spanned by the $n+1$ vertices, with the vertex labeled $i$ satisfying $t_i=1$. There are maps $\Delta_n \rightarrow \Delta_m$ induced by non-decreasing maps $\{0, \ldots, n \} \rightarrow \{0, \ldots, m \}$. The category of finite non-empty sets is called the \emph{simplex category} $\Delta$, and a simplicial object of $\msf{C}$ is a functor $\Delta^{\op} \rightarrow \msf{C}$. The corresponding functor category $\Fun(\Delta^{\op} , \msf{C})$ is abbreviated $\msf{sC}$.

There is a so-called ``Quillen equivalence'' between simplicial sets and CW complexes (a class of ``nice'' topological spaces), which informally says that one can think of simplicial sets and CW complexes as being interchangeable up to homotopy. For this reason, one often refers to simplicial sets as ``spaces'', and thinks of the adjective ``simplicial'' as synonymous to ``topological''. A {\it simplicial  group} (resp. {\it simplicial abelian group}) is a simplicial object in the category of groups (resp. abelian groups). 

The formalism of derived functors makes use, not only of  simplicial objects, but  also of appropriate generalizations of ``quasi-isomorphisms'' and ``projective resolutions''. Such notions are provided by Quillen's theory of \emph{model categories}, which is a specification of distinguished families of morphisms in $\msf{sC}$ satisfying suitable properties. The existence of a model category structure on $\msf{sC}$ is not guaranteed, but much work has gone into producing such structures on categories of interest. Often one starts with a standard model structure on the category $\msf{sSet}$ of simplicial sets, and then bootstraps from this to sets with finitary algebraic structure such as groups, abelian groups, rings, etc. We will not discuss these details; see instead \cite{DS}. Once a model category structure is in place, one constructs derived functors by a procedure analogous to the traditional calculus in derived categories, using ``projective resolutions''. (Actually in \cite{FHMR} this will all be used in a different way, using the framework of $\infty$-categories \cite{L}.)

What does this have to do with homological algebra? Recall that a chain complex is called \emph{connective} if it is supported in non-negative degrees. If $\msf{C}$ is an abelian category, then we denote by $\Ch_+(\msf{C})$ the category of connective chain complexes of objects in $\msf{C}$. The \emph{Dold-Kan correspondence} \cite[\S III.2]{GJ} gives an equivalence between $\msf{sC}$ and $\Ch_+(\msf{C})$, demonstrating that in the case of abelian categories the ``simplicial'' theory of homotopical algebra recovers the older ``chain complex'' theory of homological algebra.

Circling back to abelian groups: we use the formalism of simplicial abelian groups as the context for derived functors involving abelian groups. A simplicial abelian group $G^{\bu}$ has homotopy groups $\pi_i(G^{\bu})$; these coincide with the homotopy groups of the topological space corresponding to the underlying simplicial set of $G^{\bu}$. The ``singular simplices'' functor from topological spaces to simplicial sets promotes to a functor from topological abelian groups to simplicial abelian groups. With these preparatory remarks in place, we return to the subject of derived abelianization.

\begin{const}[Derived abelianization]\label{derab} The paper \cite{FHMR} gives several explicit descriptions of the derived abelianization $\Gamma^{\dab}$, which require a bit more language to explain. Instead, we will give a more down-to-earth model for its homotopy type. Let $\Gamma$ be a simplicial abelian group. Let $(B \Gamma, e)$ be the bar construction on $\Gamma$, viewed as a pointed space (see \cite[Chapter 16, \S 5]{May} for an explanation of the bar construction). 

Given any pointed space $(X, x)$, there is the infinite symmetric product \cite[p.282]{H} 
\[
\Sym(X,x) = \colim_n \Sym^n(X)
\]
where $\Sym^n(X) = X^n/S_n$ and the transition maps append the basepoint $x$. It is a topological abelian monoid, under concatenation. 

Then the homotopy type of the derived abelianization of $\Gamma$ is represented by the topological abelian group $\Omega \Sym(B\Gamma, e)$ where $\Omega$ is the (based) loop space functor. 

As a sanity check, note that if $\Gamma$ is discrete, then 
\[
\pi_0 (\Omega \Sym(B\Gamma, e) ) \cong \pi_1(\Sym (B\Gamma, e)) \cong H_1(\Gamma; \Z) \cong \Gamma^{\ab}.
\]
This affirms the intuition that (for discrete $\Gamma$) $\Gamma^{\dab}$ should be a space whose underlying group of connected components is $\Gamma^{\ab}$. In fact, the \emph{Dold-Thom theorem} \cite{DT} implies that---as signaled  in (\ref{h1}) above--- for all $i \geq 0$  we have:
\begin{equation}\label{eq: DT}
\pi_i(\Omega \Sym(B\Gamma, e) ) \cong \pi_{i+1}(\Sym (B\Gamma, e)) \cong H_{i+1}(\Gamma; \Z),
\end{equation}
giving us some understanding of the higher homotopy groups of $\Gamma^{\dab}$ as well. 

Under the Quillen equivalence between CW complexes and simplicial abelian groups, $\Omega \Sym (B\Gamma, e)$ may be viewed as a simplicial abelian group. Then under the Dold-Kan equivalence, it corresponds to a connective chain complex. This turns out to be a familiar object: as a chain complex, $\Gamma^{\dab}$ is quasi-isomorphic to the (reduced, shifted) homology chains $\overline{C_*(X,\Z)}[-1]$; this is a refinement of \eqref{eq: DT}. 

\end{const}

\begin{remark}[Derived abelianization of profinite groups]\label{rem: profinite abelianization} Since we are interested in Galois groups, we will want to take the derived abelianization of profinite groups. In this case it is natural to modify the derived abelianization construction to produce a profinite abelian group. If $\Gamma$ is a profinite group, then we denote by $\Gamma^{\dab}$ its profinite derived abelianization (the universal profinite simplicial abelian group to which it maps). This can be described explicitly as follows: if $\Gamma = \limit_\alpha \Gamma_{\alpha}$ is a profinite presentation of $\Gamma$, then $\Gamma^{\dab} \cong \limit_\alpha \Gamma_\alpha^{\dab}$.    
\end{remark}



 \subsection{Derived abelianization of Galois groups}
Let $\Gamma_S = \pi_1(Spec(\cO_F[1/S]))$ be the Galois group of the maximal extension $F_S/F$ unramified outside a finite set $S$ of prime ideals, equipped with its natural profinite construction. We let $\Gamma_S^{\dab}$ be as in Remark \ref{rem: profinite abelianization}. Then under a profinite version of the Dold-Kan correspondence, we have 
\begin{equation}\label{DAB}
\Gamma_S^\dab \isoarrow  \overline{C_*}(\Gamma_S,\wh{\Z})[-1].
\end{equation}
In particular,
\begin{equation}\label{pi0dab}
\pi_0(\Gamma_S^\dab) = H_1(\Gamma_S, \wh{\Z}) \cong \Gamma_S^{ab}
\end{equation}
is the classical profinite abelianization. 


\section{A derived Langlands correspondence for $\GL(1)$  -- the automorphic side }\label{derauto}
Class field theory identifies the classical abelianization of $\pi_1(\Spec(\cO_F[1/S]))$ with the class group of $\cO_F[1/S]$. We will now describe an enhancement of this story. To simplify the exposition we will focus on the case where $S$ is empty and  $F$ is totally imaginary. Otherwise, there are subleties coming from the interaction of real places with the prime $2$. These subtleties are treated in detail in \cite{FHMR}. 

\subsection{The Picard groupoid of $\Spec(\cO_F)$}\label{picg}

Recall that the class group of $\cO_F$ may be defined as the group of equivalence classes of line bundles over $\Spec \cO_F$. A more refined structure is the \emph{Picard groupoid} of $\cO_F$, which is the category whose objects are line bundles on $\cO_F$ and morphisms are isomorphisms of line bundles. This construction makes sense more generally on any scheme $X$, and we shall denote the Picard groupoid of $X$ by $\Pic(X)$. This Picard groupoid may be viewed as a simplicial set via the simplicial nerve construction; indeed a fundamental tenet of the approach to higher category theory in \cite{L} is to unify spaces and categories under the umbrella of simplicial sets. Moreover, the tensor product of line bundles induces a symmetric monoidal structure on the category $\Pic(X)$, which translates into an abelian group structure on the corresponding simplicial set. Denote by $\Pic(X)$ the Picard groupoid of $\cO_F$, viewed interchangeably as a symmetric monoidal category or as a simplicial abelian group. In these terms the class group may be described as 
\begin{equation}\label{eq: class group as pi_0}
\Cl(\cO_F)= \pi_0 \Pic(\cO_F).
\end{equation}
This is a reflection of the fact that for any scheme $X$, there is a natural isomorphism between $\rH^1(X, \G_m)$ and equivalence classes of line bundles on $X$. More generally, there is a cohomological description of the Picard groupoid. We may view $\rH^1(X, \G_m)$ as the 0th cohomology group of the cohomology complex $\rC^\bullet(X, \G_m[1])$, which is well-defined in the homotopy category of complexes. There is a truncation functor $\tau_{\leq 0}$ on the category of chain complexes, which extracts the connective cover of a complex, and it is a general fact that the Picard category of $X$ is naturally isomorphic to $\tau_{\leq 0}\rC^\bullet(X, \G_m[1])$, the connective cover of the cohomology complex $\rC^\bullet(X, \G_m[1])$. Taking the 0th cohomology group of this isomorphism recovers \eqref{eq: class group as pi_0}.

\subsection{The Picard groupoid and the id\`ele class groupoid}\label{idele}
If $F$ is the function field of a curve $X$ over a finite field, then Weil's construction identifies the id\`ele class group $I_{F,K_{max}}$, where $K_{max}$ is the product over all places $v$ of $F$ of the maximal compact subgroups of $F_v^\times$, with the groupoid of line bundles on $X$. A similar construction can be applied when $F$ is a number field, which involves metrics at archimedean places. However, these can be ignored under our simplifying assumption  that $F$ is totally imaginary. So we can get away (under  this assumption) with just considering the Picard groupoid $\Pic(\Spec(\cO_F))$; we then get a natural homotopy equivalence 
$$\Pic(\Spec(\cO_F)) \simeq I_{F,K_{max}}$$
We let $K = K_{max}$ in what follows. We then have a cohomological description of the idele class group, as 
\begin{equation}\label{picidele}
\tau_{\leq 0}\rC^\bullet(X, \G_m[1])	 \simeq \Pic(\Spec(\cO_F)) \simeq I_{F,K}.
\end{equation}


\subsection{The profinite completion of the Picard groupoid and flat cohomology}
The approach to the derived abelianization outlined above is by way of the $p$-adic derived deformation ring.  To compare this with the Picard groupoid, we introduce the $p$-adic completion of the latter.  We consider the Kummer exact sequence
\begin{equation}\label{kummer}
1 \ra \mu_{n} \ra \G_m \overset{n}\to \G_m \ra 1
\end{equation}
which we can write as an isomorphism in the derived category of sheaves on the flat topology of $\Spec(\cO_F)$:
\begin{equation}\label{kummer2}
[\G_m \overset{n}\to \G_m] \simeq \mu_{n}[1];  ~~~(\G_m)\hat  ~:=  \varprojlim_m \G_m/n  \simeq \varprojlim_m \mu_{n}[1].
\end{equation}
We denote by $C^*(\cO_F, \bullet)$ the cohomology complex of $\Spec(\cO_F)$ of the fppf sheaf $\bullet$. The profinite completion of $C^*(X, \G_m)$ is therefore 
\begin{equation}\label{pcomplete}
C^*(\cO_F, \G_m)\hat{} := \varprojlim_n C^*(\cO_F, \G_m)/n \isoarrow C^*(\cO_F, \varprojlim_n \G_m/n) \cong C^*(\cO_F, \mu[1])
\end{equation}
where $\mu := \varprojlim_n \mu_n$ is the Tate module of roots of unity. 

Combining \eqref{pcomplete} with \eqref{kummer2}, we identify the profinite completion of $\rC^*(\cO_F, \G_m)$ with the flat cohomology complex:
\begin{equation}\label{kummer3}
\rC^*(\cO_F, \G_m)\hat{}  ~~~\simeq ~~~ C^*(\cO_F,\mu[1])
\end{equation}
There is a generalization of profinite completion to simplicial sets and simplicial abelian groups. The considerations above then yield 
\begin{equation}\label{kummer4}
\Pic(\cO_F)\hat{}  \isoarrow  \tau_{\leq 0} C^*(\cO_F, \mu[2]).
\end{equation}

\section{Derived class field theory via Poitou-Tate duality}\label{PTsection}

\subsection{Derived Poitou-Tate duality}

Recall that an oriented manifold $M$ enjoys Poincar\'{e} duality, which can be formulated as an isomorphism 
\[
H_i(M,\Z) \cong H^{n-i}_c(M, \Z) 
\]
where $n = \dim M$. In fact, this can be promoted to an isomorphism of complexes (in a suitable localization of the category of complexes)
\[
C_i(M, \Z) \cong C^{n-i}_c(M, \Z).
\]
There is an analogy between number fields and 3-manifolds, under which Poincar\'{e} duality is analogous to the so-called Poitou-Tate duality. The latter is a bit subtle, but the upshot is that for a finite fppf sheaf $\cF$ with Cartier dual $\cF^D := \cHom(\cF, \G_m)$, there are isomorphisms 
\[
H_i(\cO_F, \cF) \cong H^{3-i}_c(\cO_F, \cF^D).
\]
Here the compactly supported cohomology groups $ H^{3-i}_c(\cO_F, \cF^D)$ are a bit more involved to define in general -- they are defined formally as cones of restriction maps to Archimedean places, where one also has to replace cohomology by Tate cohomology. However, under our simplifying assumptions we will have $H^{3-i}_c(\cO_F, \cF^D) = H^{3-i}(\cO_F, \cF^D)$, and we can ignore the compact support condition entirely. It turns out that in this case, similarly to the case of manifolds, one can promote Poitou-Tate duality to an isomorphism 
\[
C_*(\cO_F, \cF) \cong C^*(\cO_F, \cF^D[3]). 
\]
This promotion is quite formal but we prefer to leave the details to \cite{FHMR}. Applying this with $\cF = \mu_n$ and taking limits in $n$, the upshot is a natural isomorphism 
\[
C_*(\cO_F, \wh{\Z}) \cong C^*(\cO_F,  \mu[3]).
\]
Combining this with \eqref{kummer4} gives the following:
\begin{equation}\label{PT3}
\Pic(\cO_F)\hat{} \isoarrow \tau_{\leq 0}C^*_{flat}(\cO_F,  \mu[2]) \isoarrow \tau_{\leq 0} C_*(\pi_1(\cO_F), \wh{\Z}[-1])
\end{equation}

\subsection{Derived class field theory}

Now we put together the isomorphism \eqref{PT3} with \eqref{DAB} and \eqref{kummer3} to find the following dual description of the derived abelianization of Galois groups: 

\begin{thm}\label{mainthm_noS}  Suppose that $F$ is a totally imaginary number field. Then there is a natural isomorphism
$$\left( \pi_1(\Spec(\cO_F))^\dab \right)\hat{} \simeq  (I_{F,K})\hat{}$$
where $K = \prod_v \cO_{F,v}^\times$ is the maximal compact subgroup of the finite id\`eles, such that taking $\pi_0$ recovers the usual isomorphism of class field theory 
\[
\pi_1(\Spec(\cO_F))^{\ab} \cong \pi_0 (I_{F,K}).
\]
\end{thm}

\subsection{Allowing for ramification}

The version of derived class field theory explained above recovers the id\`ele quotient $\pi_0(I_{F,K})$ when $K$ is the product of the unit groups $\cO_v$ at all primes $v$.  For more general $K$, one can interpret $\pi_0(I_{F,K})$ as the group of equivalence classes of invertible $\cO_F$-modules $M$ with a given {\it level $K$ structure}. This is defined as follows.  If $K \subset \prod_{v \in S} \cO_v$
is the principal congruence subgroup $K_J$ of level $J$, where $J \subset \cO_F$ is an ideal supported at $S$, then a level $K$ structure is a trivialization
$$\iota:  \cO_F/J \isoarrow M/JM.$$
In general, any open subgroup $K \subset \prod_{v \in S} \cO_v$ contains some $K_J$, and a level $K$ structure is an equivalence class of level $K_J$-structures for the action of $K/K_J$.  There is a (relatively formal) way to enhance Theorem \ref{mainthm_noS} to encompass such level structures, using \emph{relative derived abelianization}, which will be explained in \cite{FHMR}. 

\section{Functoriality}  

\subsection{Behavior under finite field extensions}  
Let $E/F$ be a finite extension of degree $n$.    Fix a finite set $S$ of non-archimedean places of $F$ and let $\Gamma_{F,S} =  \pi_1(\cO_F)$,
$\Gamma_{E,S} =  \pi_1(\cO_E)$.  Let $K^S_F = \prod_{v \notin S} \cO_v \subset \A_F^{\times}$, $K^S_E \subset \A_E^{\times}$, be the corresponding subgroups
The following commutative diagrams are the derived versions of the familiar functorialities of class field theory; the profinite completions have been omitted.

\begin{equation}\label{norm}
\xymatrix{
I_{K^S_E}  \ar[d]^{N_{E/F}} \ar[r]   &\Gamma_{E,S}^\dab \ar[d]^{incl} \\
I_{K^S_F}  \ar[r]   &\Gamma_{F,S}^\dab 
}
\end{equation}

\begin{equation}\label{transfer}
\xymatrix{
I_{K^S_E}  \ar[r]   &\Gamma_{E,S}^\dab  \\
I_{K^S_F}\ar[u]_{incl}  \ar[r]   &\Gamma_{F,S}^\dab \ar[u]_{tr}
}
\end{equation}

The horizontal maps are the derived class field theory equivalences discussed above. 

The vertical maps labelled $incl$ are the natural inclusions, which
go in opposite directions for the id\`ele classes and Galois groups (inclusion followed by abelianization, in the latter case).  The map $N_{E/F}$ is the norm map.  As usual, the only map that requires explanation is the right-hand vertical map in  \eqref{transfer}.  This is the transfer, which on homotopy groups is given by the corestriction on group cohomology. 

\subsection{Iwasawa theory}

Now suppose $F = F_0 \subset F_1 \subset \dots \subset F_n \subset \dots$ is a tower of extensions of $F$, with 
$\Gal(F_n/F) \isoarrow \ZZ/p^n\ZZ$, so that $F_\infty = \cup_n F_n$ is a {\it $\ZZ_p$-extension} of $F$.  This is the setting of {\it Iwasawa theory} to which
so much of John Coates's work was devoted.  Write $\Gamma = \Gal(F_\infty/F)$. 

In classical Iwasawa theory, one considers the limit and colimit of the classical abelianizations $\Gamma_{F_n,S,p}^{\ab}$ in this tower
\begin{equation}\label{Iwa}
\Gamma_{F_\infty,S,p}^{\ab} := \varprojlim_n \Gamma_{F_n,S,p}^{\ab} ;   ~~~A_{F_\infty,S,p}^{\ab} = \varinjlim_n \Gamma_{F_n,S,p}^{\ab} .
\end{equation}
Here the subscript $p$ denotes $p$-adic completion. One views the resulting objects in \eqref{Iwa} as compact and discrete modules, respectively, over the Iwasawa algebra $\ZZ_p[[\Gamma]]$, where the double brackets denote the completion of the group algebra with respect to the inverse limit topology. 

There is a remarkable analogy between the structure of $A_{F_\infty,S,p}^{\ab}$ as a module over the Iwasawa algebra $\ZZ_p[[\Gamma]]$ and the geometric cohomology of a smooth projective curve $X/\F_p$ as a module over the Frobenius, which is at the heart of Iwasawa theory. In the setting of the smooth projective curve, recent work in the Geometric Langlands program \cite{AGKRRV} proposes a refined Langlands correspondence for $X$, using the realization of the derived space of local systems on $X$ as the Frobenius-fixed points of the derived space of local systems on $X_{\ol{\F}_p}$.  (See also \cite{Gaits}.) The incorporation of derived structure is crucial, and one can speculate that a similar story exists over number fields, with the Iwasawa algebra playing the role of the Frobenius. This speculation is currently being investigated in work-in-progress of the authors.


\end{document}